\documentclass{article}
\usepackage{amsfonts}
\usepackage{amsmath}

\input{tcilatex}
\begin{document}

\title{The $m$-th root Finsler geometry of the Bogoslovsky-Goenner metric}
\author{Mircea Neagu}
\date{}
\maketitle

\begin{abstract}
In this paper we present the $m$-th root Finsler geometries of the three and
four dimensional Bogoslovsky-Goenner metrics (good Finslerian anisotropic
models in Special Relativity), in the sense of their Cartan torsion and
curvature distinguished tensors or vertical Einstein-like equations.
\end{abstract}



\section{Introduction}

The physical studies of Asanov \cite{Asanov}, Bogoslovsky \cite{Bogoslovsky}
or Minguzzi \cite{Minguzzi} emphasize the importance of the Finsler geometry
in the theory of space-time structure, gravitation and electromagnetism. For
this reason, one emphasizes the important role played by the Finsler-Asanov
metric (see Miron et al., \cite{Mir-Hr-Shi-Sab}, pp. 54)

\begin{equation*}
F:TM\rightarrow \mathbb{R},\quad F(x,y)=\left( a^{1}a^{2}...a^{n}\right) ^{%
\frac{1}{n}},
\end{equation*}%
where $a^{\alpha }$ $(\alpha =1,2,...,n)$ are linearly independent 1-forms.
The above Finsler metric was initially considered, in diverse particular
forms, by Riemann and Berwald-Mo\'{o}r (see Minguzzi \cite{Minguzzi} or
Bogoslovsky and Goenner \cite{Bog-G} and references therein). Moreover,
considering that $a^{\alpha }=a_{\beta }^{\alpha }y^{\beta }$, where $%
a_{\beta }^{\alpha }\in \mathbb{R}$, the above Finsler-Asanov metric becomes
a Minkowski metric of $n$-th root type. These kind of metrics were
intensively studied by Shimada \cite{Shimada}.

In such an anisotropic physical context, Bogoslovski and Goenner introduced
the locally Minkowski metric of a flat space-time with entirely broken
isotropy, which is given by

\begin{enumerate}
\item for $n=3$:%
\begin{equation}
L(y)=\sqrt[3]{(y^{1}-y^{2}-y^{3})(y^{1}-y^{2}+y^{3})(y^{1}+y^{2}-y^{3})};
\label{n=3}
\end{equation}

\item for $n=4$:%
\begin{eqnarray}
L(y) &=&\sqrt[4]{(y^{1}-y^{2}-y^{3}-y^{4})(y^{1}-y^{2}+y^{3}+y^{4})\cdot }
\label{n=4} \\
&&\overline{\cdot (y^{1}+y^{2}-y^{3}+y^{4})(y^{1}+y^{2}+y^{3}-y^{4})}. 
\notag
\end{eqnarray}
\end{enumerate}

\section{The Finsler geometry of a locally Minkovski space}

Let us consider a locally Minkowski space $(M^{n},L=L(y))$, where $%
L:TM\rightarrow \mathbb{R}$ is a Finsler metric depending only on
directional variables $y^{i}$, where $i=1,2,...,n$. Note that, in this
Section, the Latin letters $i,j,k,...$ run from $1$ to $n$, and the Einstein
convention of summation is adopted all over this work. It follows that the
fundamental metrical d-tensor of the Minkowski space is%
\begin{equation*}
g_{ij}(y)=\frac{1}{2}\frac{\partial ^{2}L^{2}}{\partial y^{i}\partial y^{j}},
\end{equation*}%
whose inverse d-tensor is given by $g^{jk}(y)$. Taking into account that the
Finsler function $L$ is a locally Minkowski metric, we deduce that the
Euler-Lagrange equations of $L$ produce the canonical nonlinear connection $%
N_{j}^{i}=0$, where%
\begin{equation*}
N_{j}^{i}=\frac{\partial G^{i}}{\partial y^{j}}=\frac{\partial }{\partial
y^{j}}\left[ g^{iu}\left( \frac{\partial ^{2}L^{2}}{\partial x^{v}\partial
y^{u}}y^{v}-\frac{\partial L^{2}}{\partial x^{u}}+\frac{\partial ^{2}L^{2}}{%
\partial t\partial y^{u}}\right) \right] .
\end{equation*}%
As a consequence, in the Finsler geometrical study of the Minkowski metric
are important only the vertical geometrical objects like (see \cite%
{Mir-Hr-Shi-Sab}):

\begin{enumerate}
\item the Cartan d-torsion:%
\begin{equation*}
C_{jkm}=\frac{1}{2}\frac{\partial g_{jk}}{\partial y^{m}}=\frac{1}{4}\frac{%
\partial ^{3}L^{2}}{\partial y^{j}\partial y^{k}\partial y^{m}}\quad \text{%
(covariant form),}
\end{equation*}%
\begin{equation*}
C_{jk}^{i}=g^{im}C_{jkm}=\frac{g^{im}}{2}\frac{\partial g_{jk}}{\partial
y^{m}}\quad \text{(contravariant form);}
\end{equation*}

\item the Cartan d-curvature:%
\begin{equation*}
S_{ijk}^{l}=\frac{\partial C_{ij}^{l}}{\partial y^{k}}-\frac{\partial
C_{ik}^{l}}{\partial y^{j}}+C_{ij}^{u}C_{uk}^{l}-C_{ik}^{u}C_{uj}^{l}\quad 
\text{(contravariant form),}
\end{equation*}%
\begin{equation*}
S_{imjk}=g_{ml}S_{ijk}^{l}=g^{uv}\left( C_{ujm}C_{vik}-C_{ukm}C_{vij}\right)
\quad \text{(covariant form).}
\end{equation*}

\item the vertical Einstein-like equations for $n>2$:%
\begin{equation*}
S_{ij}-\frac{S}{2}g_{ij}=\text{$\QTR{sc}{k}$}\widetilde{T}_{ij},
\end{equation*}%
where

\begin{itemize}
\item $S_{ij}=S_{ijm}^{m}$ is the vertical Ricci d-tensor;

\item $S=g^{uv}S_{uv}$ is the scalar curvature;

\item $\widetilde{T}_{ij}$ are the new vertical components of the
non-isotropic stress-energy d-tensor of matter $\mathbb{T}$;

\item $\QTR{sc}{k}$ is the Einstein constant.
\end{itemize}
\end{enumerate}

\section{The $3$-rd root Finsler geometry of the three dimensional
Bogoslovsky-Goenner metric}

In this Section we have $M=\mathbb{R}^{3}$, that is $n=3$, and the Latin
indices $i,j,k,...$ run from $1$ to $3$. Let us consider the notations $%
S_{\alpha }=(y^{1})^{\alpha }+(y^{2})^{\alpha }+(y^{3})^{\alpha }$, where $%
\alpha \in \mathbb{Z}$, $P_{3}=y^{1}y^{2}y^{3}$ and 
\begin{equation*}
A=(y^{1}-y^{2}-y^{3})(y^{1}-y^{2}+y^{3})(y^{1}+y^{2}-y^{3}).
\end{equation*}%
In such a context, by direct computations, the three dimensional
Bogoslovsky-Goenner metric (\ref{n=3}) takes the 3-rd root metric form%
\begin{equation*}
L=\sqrt[3]{A}=\sqrt[3]{%
(y^{1})^{3}+(y^{2})^{3}+(y^{3})^{3}-y^{1}(y^{2})^{2}-y^{1}(y^{3})^{2}-y^{2}(y^{1})^{2}-y^{2}(y^{3})^{2}-%
}
\end{equation*}%
\begin{equation*}
\overline{-y^{3}(y^{1})^{2}-y^{3}(y^{2})^{2}+2y^{1}y^{2}y^{3}}=\sqrt[3]{%
2S_{3}-S_{1}S_{2}+2P_{3}}.
\end{equation*}

Working on the domain in which $2S_{3}-S_{1}S_{2}+2P_{3}\neq 0$, the
fundamental metrical d-tensor produced by the Bogoslovsky-Goenner metric of
order three (\ref{n=3}) is given by%
\begin{equation}
g_{ij}=\frac{1}{2}\frac{\partial ^{2}L^{2}}{\partial y^{i}\partial y^{j}}=%
\frac{1}{2}\frac{\partial ^{2}A^{2/3}}{\partial y^{i}\partial y^{j}}=\frac{%
A^{-1/3}}{3}A_{ij}-\frac{A^{-4/3}}{9}A_{i}A_{j},  \label{g_ij3}
\end{equation}%
where%
\begin{equation*}
A_{i}=\frac{\partial A}{\partial y^{i}}=6(y^{i})^{2}-S_{2}-2y^{i}S_{1}+2%
\frac{P_{3}}{y^{i}},
\end{equation*}%
\begin{equation*}
A_{ij}=\frac{\partial ^{2}A}{\partial y^{i}\partial y^{j}}=-2y^{i}-2y^{j}+2%
\frac{P_{3}}{y^{i}y^{j}}+\left[ 12y^{i}-2S_{1}-2\frac{P_{3}}{(y^{i})^{2}}%
\right] \delta _{ij}.
\end{equation*}

Putting the coefficients $A_{ij}$ into a matrix, we get the matrix%
\begin{equation*}
\left( A_{ij}\right) _{i,j=\overline{1,3}}=\left( 
\begin{array}{ccc}
6y^{1}-2y^{2}-2y^{3} & -2y^{1}-2y^{2}+2y^{3} & -2y^{1}+2y^{2}-2y^{3} \\ 
-2y^{1}-2y^{2}+2y^{3} & -2y^{1}+6y^{2}-2y^{3} & 2y^{1}-2y^{2}-2y^{3} \\ 
-2y^{1}+2y^{2}-2y^{3} & 2y^{1}-2y^{2}-2y^{3} & -2y^{1}-2y^{2}+6y^{3}%
\end{array}%
\right) ,
\end{equation*}%
whose determinant is given by $\det (A_{ij})=-8D$, where $%
D=-4(y^{1})^{3}-4(y^{2})^{3}-4(y^{3})^{3}+4(y^{1})^{2}y^{2}+4(y^{1})^{2}y^{3}+4y^{1}(y^{2})^{2}+4y^{1}(y^{3})^{2}+4(y^{2})^{2}y^{3}+4y^{2}(y^{3})^{2}-8y^{1}y^{2}y^{3}=-8S_{3}+4S_{1}S_{2}-8P_{3} 
$. If we have $D\neq 0$, then the inverse matrix of $(A_{ij})$ is the matrix 
$(A^{jk})_{j,k=\overline{1,3}}=D^{-1}\cdot A^{\ast }$, where $A^{\ast }=$%
\begin{equation*}
={\small {\left( 
\begin{array}{ccc}
2(y^{2})^{2}-4y^{2}y^{3}+2(y^{3})^{2} & S_{2}-2y^{1}y^{3}-2y^{2}y^{3} & 
S_{2}-2y^{1}y^{2}-2y^{2}y^{3} \\ 
S_{2}-2y^{1}y^{3}-2y^{2}y^{3} & 2(y^{1})^{2}-4y^{1}y^{3}+2(y^{3})^{2} & 
S_{2}-2y^{1}y^{2}-2y^{1}y^{3} \\ 
S_{2}-2y^{1}y^{2}-2y^{2}y^{3} & S_{2}-2y^{1}y^{2}-2y^{1}y^{3} & 
2(y^{1})^{2}-4y^{1}y^{2}+2(y^{2})^{2}%
\end{array}%
\right) .}}
\end{equation*}%
As a general formula, we have%
\begin{equation*}
A^{jk}=\frac{1}{D}\left[ S_{2}-2(y^{j}+y^{k})\frac{P_{3}}{y^{j}y^{k}}%
-(y^{j})^{2}\delta _{jk}\right] .
\end{equation*}

Using the above geometrical entities, we deduce that the inverse d-tensor of
the fundamental metrical d-tensor (\ref{g_ij3}) has the form%
\begin{equation}
g^{jk}=3A^{1/3}A^{jk}+\frac{A^{-2/3}}{1-\frac{A^{-1}}{3}A^{uv}A_{u}A_{v}}%
A^{j}A^{k},  \label{g^jk3}
\end{equation}%
where $A^{j}=A^{jw}A_{w}$. Moreover, by direct computations, the covariant
Cartan d-torsion produced by the three dimensional Bogoslovsky-Goenner
metric (\ref{n=3}) is given by%
\begin{equation*}
C_{jkm}=\frac{A^{-1/3}}{6}A_{jkm}-\frac{A^{-4/3}}{18}%
(A_{jk}A_{m}+A_{km}A_{j}+A_{mj}A_{k})-\frac{A^{-7/3}}{18}A_{j}A_{k}A_{m},
\end{equation*}%
where%
\begin{equation*}
A_{jkm}=\frac{\partial A_{jk}}{\partial y^{m}}=\frac{\partial ^{3}A}{%
\partial y^{j}\partial y^{k}\partial y^{m}}=\left\{ 
\begin{array}{ll}
6, & \text{if }j=k=m \\ 
2, & \text{if }j\neq k\neq m\neq j \\ 
-2, & \text{otherwise.}%
\end{array}%
\right.
\end{equation*}%
In other words, if we denote by $A_{(m)}=\left( A_{jkm}\right) _{j,k=%
\overline{1,3}}$, where $m\in \{1,2,3\}$, we get%
\begin{equation*}
A_{(1)}=\left( 
\begin{array}{ccc}
6 & -2 & -2 \\ 
-2 & -2 & 2 \\ 
-2 & 2 & -2%
\end{array}%
\right) ,\quad A_{(2)}=\left( 
\begin{array}{ccc}
-2 & -2 & 2 \\ 
-2 & 6 & -2 \\ 
2 & -2 & -2%
\end{array}%
\right) ,
\end{equation*}%
\begin{equation*}
A_{(3)}=\left( 
\begin{array}{ccc}
-2 & 2 & -2 \\ 
2 & -2 & -2 \\ 
-2 & -2 & 6%
\end{array}%
\right) .
\end{equation*}

\section{The $4$-th root Finsler geometry of the four dimensional
Bogoslovsky-Goenner metric}

In this Section we have $M=\mathbb{R}^{4}$, that is $n=4$, and the Latin
indices $i,j,k,...$ run from $1$ to $4$. Let us consider the notations $%
S_{\alpha }=(y^{1})^{\alpha }+(y^{2})^{\alpha }+(y^{3})^{\alpha
}+(y^{4})^{\alpha }$, where $\alpha \in \mathbb{Z}$, $%
P_{4}=y^{1}y^{2}y^{3}y^{4}$ and%
\begin{equation*}
A=(y^{1}-y^{2}-y^{3}-y^{4})(y^{1}-y^{2}+y^{3}+y^{4})(y^{1}+y^{2}-y^{3}+y^{4})(y^{1}+y^{2}+y^{3}-y^{4}).
\end{equation*}%
In such a context, by direct computations, the four dimensional
Bogoslovsky-Goenner metric (\ref{n=4}) takes the 4-th root metric form%
\begin{equation*}
L=\sqrt[4]{A}=\sqrt[4]{(y^{1})^{4}+(y^{2})^{4}+(y^{3})^{4}+(y^{4})^{4}-2%
\left( y^{1}\right) ^{2}(y^{2})^{2}-2\left( y^{1}\right) ^{2}(y^{3})^{2}-}
\end{equation*}%
\begin{equation*}
\overline{-2\left( y^{1}\right) ^{2}(y^{4})^{2}-2\left( y^{2}\right)
^{2}(y^{3})^{2}-2\left( y^{2}\right) ^{2}(y^{4})^{2}-2\left( y^{3}\right)
^{2}(y^{4})^{2}-8y^{1}y^{2}y^{3}y^{4}}=
\end{equation*}%
\begin{equation*}
=\sqrt[4]{2S_{4}-S_{2}^{2}-8P_{4}}.
\end{equation*}

Working on the domain in which $2S_{4}-S_{2}^{2}-8P_{4}>0$, the fundamental
metrical d-tensor produced by the Bogoslovsky-Goenner metric of order four (%
\ref{n=4}) is given by%
\begin{equation}
g_{ij}=\frac{1}{2}\frac{\partial ^{2}L^{2}}{\partial y^{i}\partial y^{j}}=%
\frac{1}{2}\frac{\partial ^{2}A^{1/2}}{\partial y^{i}\partial y^{j}}=\frac{%
A^{-1/2}}{4}A_{ij}-\frac{A^{-3/2}}{8}A_{i}A_{j},  \label{g_ij4}
\end{equation}%
where%
\begin{equation*}
A_{i}=\frac{\partial A}{\partial y^{i}}=4(y^{i})^{3}-4y^{i}\left[
S_{2}-\left( y^{i}\right) ^{2}\right] -8\frac{P_{4}}{y^{i}},
\end{equation*}%
\begin{equation*}
A_{ij}=\frac{\partial ^{2}A}{\partial y^{i}\partial y^{j}}=-8y^{i}y^{j}-8%
\frac{P_{4}}{y^{i}y^{j}}+\left[ 24\left( y^{i}\right) ^{2}-4S_{2}+8\frac{%
P_{4}}{(y^{i})^{2}}\right] \delta _{ij}.
\end{equation*}

Putting the coefficients $A_{ij}$ into a matrix, we get the matrix%
\begin{equation*}
\left( A_{ij}\right) _{i,j=\overline{1,4}}=\left( 
\begin{array}{cccc}
p & a & b & c \\ 
a & q & c & b \\ 
b & c & r & a \\ 
c & b & a & s%
\end{array}%
\right) ,
\end{equation*}%
where%
\begin{equation*}
p=12(y^{1})^{2}-4(y^{2})^{2}-4(y^{3})^{2}-4(y^{4})^{2},
\end{equation*}%
\begin{equation*}
q=-4(y^{1})^{2}+12(y^{2})^{2}-4(y^{3})^{2}-4(y^{4})^{2},
\end{equation*}%
\begin{equation*}
r=-4(y^{1})^{2}-4(y^{2})^{2}+12(y^{3})^{2}-4(y^{4})^{2},
\end{equation*}%
\begin{equation*}
s=-4(y^{1})^{2}-4(y^{2})^{2}-4(y^{3})^{2}+12(y^{4})^{2},
\end{equation*}%
\begin{equation*}
a=-8y^{1}y^{2}-8y^{3}y^{4},
\end{equation*}%
\begin{equation*}
b=-8y^{1}y^{3}-8y^{2}y^{4},
\end{equation*}%
\begin{equation*}
c=-8y^{1}y^{4}-8y^{2}y^{3},
\end{equation*}%
whose determinant is given by $D=\det
(A_{ij})=a^{4}-2a^{2}c^{2}-2b^{2}c^{2}-2a^{2}b^{2}+b^{4}+c^{4}-a^{2}pq-b^{2}pr-a^{2}rs-b^{2}qs-c^{2}ps-c^{2}qr+2abcp+2abcq+2abcr+2abcs+pqrs 
$. If we have $D\neq 0$, then the inverse matrix of $(A_{ij})$ is the matrix 
$(A^{jk})_{j,k=\overline{1,4}}=D^{-1}\cdot A^{\ast }$, where $A^{\ast }=$%
\begin{equation*}
=\left( 
\begin{array}{cc}
-qa^{2}+2abc-rb^{2}-sc^{2}+qrs & a^{3}-ac^{2}-ab^{2}+bcr+bcs-ars \\ 
a^{3}-ac^{2}-ab^{2}+bcr+bcs-ars & -pa^{2}+2abc-sb^{2}-rc^{2}+prs \\ 
b^{3}-bc^{2}-a^{2}b+acq+acs-bqs & c^{3}-b^{2}c-a^{2}c+abp+abs-cps \\ 
c^{3}-b^{2}c-a^{2}c+abq+abr-cqr & b^{3}-bc^{2}-a^{2}b+acp+acr-bpr%
\end{array}%
\right.
\end{equation*}%
\begin{equation*}
\left. 
\begin{array}{cc}
b^{3}-bc^{2}-a^{2}b+acq+acs-bqs & c^{3}-b^{2}c-a^{2}c+abq+abr-cqr \\ 
c^{3}-b^{2}c-a^{2}c+abp+abs-cps & b^{3}-bc^{2}-a^{2}b+acp+acr-bpr \\ 
-sa^{2}+2abc-pb^{2}-qc^{2}+pqs & a^{3}-ac^{2}-ab^{2}+bcp+bcq-apq \\ 
a^{3}-ac^{2}-ab^{2}+bcp+bcq-apq & -ra^{2}+2abc-qb^{2}-pc^{2}+pqr%
\end{array}%
\right) .
\end{equation*}

Using the above geometrical entities, we deduce that the inverse d-tensor of
the fundamental metrical d-tensor (\ref{g_ij4}) has the form%
\begin{equation}
g^{jk}=4A^{1/2}A^{jk}+\frac{A^{-1/2}}{2-A^{-1}A^{uv}A_{u}A_{v}}A^{j}A^{k},
\label{g^jk4}
\end{equation}%
where $A^{j}=A^{jw}A_{w}$. Moreover, by direct computations, the covariant
Cartan d-torsion produced by the four dimensional Bogoslovsky-Goenner metric
(\ref{n=4}) is given by%
\begin{equation*}
C_{jkm}=\frac{A^{-1/2}}{8}A_{jkm}-\frac{A^{-3/2}}{16}%
(A_{jk}A_{m}+A_{km}A_{j}+A_{mj}A_{k})+\frac{3A^{-5/2}}{32}A_{j}A_{k}A_{m},
\end{equation*}%
where%
\begin{equation*}
A_{jkm}=\frac{\partial A_{jk}}{\partial y^{m}}=\frac{\partial ^{3}A}{%
\partial y^{j}\partial y^{k}\partial y^{m}}.
\end{equation*}%
If we denote by $A_{(m)}=\left( A_{jkm}\right) _{j,k=\overline{1,4}}$, where 
$m\in \{1,2,3,4\}$, we get%
\begin{equation*}
A_{(1)}=\left( 
\begin{array}{cccc}
24y^{1} & -8y^{2} & -8y^{3} & -8y^{4} \\ 
-8y^{2} & -8y^{1} & -8y^{4} & -8y^{3} \\ 
-8y^{3} & -8y^{4} & -8y^{1} & -8y^{2} \\ 
-8y^{4} & -8y^{3} & -8y^{2} & -8y^{1}%
\end{array}%
\right) ,
\end{equation*}%
\begin{equation*}
A_{(2)}=\left( 
\begin{array}{cccc}
-8y^{2} & -8y^{1} & -8y^{4} & -8y^{3} \\ 
-8y^{1} & 24y^{2} & -8y^{3} & -8y^{4} \\ 
-8y^{4} & -8y^{3} & -8y^{2} & -8y^{1} \\ 
-8y^{3} & -8y^{4} & -8y^{1} & -8y^{2}%
\end{array}%
\right) ,
\end{equation*}%
\begin{equation*}
A_{(3)}=\left( 
\begin{array}{cccc}
-8y^{3} & -8y^{4} & -8y^{1} & -8y^{2} \\ 
-8y^{4} & -8y^{3} & -8y^{2} & -8y^{1} \\ 
-8y^{1} & -8y^{2} & 24y^{3} & -8y^{4} \\ 
-8y^{2} & -8y^{1} & -8y^{4} & -8y^{3}%
\end{array}%
\right) ,
\end{equation*}%
\begin{equation*}
A_{(4)}=\left( 
\begin{array}{cccc}
-8y^{4} & -8y^{3} & -8y^{2} & -8y^{1} \\ 
-8y^{3} & -8y^{4} & -8y^{1} & -8y^{2} \\ 
-8y^{2} & -8y^{1} & -8y^{4} & -8y^{3} \\ 
-8y^{1} & -8y^{2} & -8y^{3} & 24y^{4}%
\end{array}%
\right) .
\end{equation*}%

Mircea Neagu

Department of Mathematics and Informatics

Transilvania University of Bra\c{s}ov

Blvd. Iuliu Maniu, No. 50, Bra\c{s}ov 500091, Romania

email: \textit{mircea.neagu@unitbv.ro}

\end{document}